\documentclass{amsart}

\usepackage{amsfonts, amsmath, amssymb, amscd, amsthm}
\usepackage{graphicx}

\newtheorem{thm}{Theorem}[section]
\newtheorem{cor}[thm]{Corollary}
\newtheorem{lem}[thm]{Lemma}

\theoremstyle{definition}
\newtheorem{defn}[thm]{Definition}
\theoremstyle{remark}

\newtheorem{ex}[thm]{Example}

\newcommand{\Hl }{\{ H_\lambda \} _{ \lambda \in \Lambda  }}
\newcommand{\Ul }{\{ U_\lambda \} _{ \lambda \in \Lambda  }}
\newcommand{\Am }{\{ A_\mu  \} _{ \mu \in {\rm M} }}
\newcommand{\Bn }{\{ B_\nu \} _{ \nu\in {\rm N}  }}
\newcommand{\G }{\Gamma }
\newcommand{\Ga }{\Gamma (G, \mathcal A)}

\newcommand{\GG }{\Gamma (G, X\cup \mathcal H)}

\renewcommand{\t }{{\rm\bf top}}
\renewcommand{\b }{{\rm\bf bot}}

\newfont{\eufm}{eufm10}

\begin{document}

\author{D.V. Osin}
\thanks{This work has been supported by the RFBR Grant $\sharp $
02-01-00892}

\address{Department of Mathematics, 1326 Stevenson Center,
Vanderbilt University, Nashville  TN 37240-0001, USA}
\email{denis.ossine@math.vanderbilt.edu \\
denis.osin@mtu-net.ru}


\title{Relative Dehn functions of amalgamated products and HNN--extensions.}

\begin{abstract}
We obtain an upper bound for relative Dehn functions of
amalgamated products and HNN--extensions with respect to certain
collections of subgroups. Our main results generalize the
combination theorems for relatively hyperbolic groups proved by
Dahmani.
\end{abstract}

\maketitle

\section{Introduction}

Relative Dehn functions were introduced in \cite{RHG} (a
particular case was also considered in \cite{BC}) in order to
obtain an isoperimetric characterization of relatively hyperbolic
groups. In this paper we prove two theorems concerning relative
Dehn functions of amalgamated products and HNN--extensions. In
particular, our results generalize the Combination Theorem for
relatively hyperbolic groups proved by Dahmani \cite{Dah}. Results
of similar type for ordinary Dehn functions can be found in
\cite{BF,BC,BC1,KM,MO}.

We begin with definitions. Let $G$ be a group, $\Hl $ a collection
of subgroups of $G$, $X$ a subset of $G$. We say that $X$ is a
{\it relative generating set of $G$ with respect to $\Hl $} if $G$
is generated by $X$ together with the union of all $H_\lambda $.
(We always assume $X$ to be symmetrized, i.e. $X^{-1}=X$.) In this
situation the group $G$ can be regarded as a quotient group of the
free product
\begin{equation}
F=\left( \ast _{\lambda\in \Lambda } H_\lambda  \right) \ast F(X),
\label{F}
\end{equation}
where $F(X)$ is the free group with the basis $X$. Let $N$ denote
the kernel of the natural homomorphism $F\to G$. If $N$ is a
normal closure of a subset $\mathcal R\subseteq N$ in the group
$F$, we say that $G$ has {\it relative presentation}
\begin{equation}\label{G}
\langle X,\; H_\lambda, \lambda\in \Lambda \; |\; \mathcal R
\rangle .
\end{equation}
If $\sharp\, X<\infty $ and $\sharp\, \mathcal R<\infty $, the
relative presentation (\ref{G}) is said to be {\it finite} and the
group $G$ is said to be {\it finitely presented relative to the
collection of subgroups $\Hl $.}

Let
\begin{equation}\label{H}
\mathcal H=\bigsqcup\limits_{\lambda\in \Lambda} (H_\lambda
\setminus \{ 1\} ) .
\end{equation}
Given a word $W$ in the alphabet $X\cup \mathcal H$ such that $W$
represents $1$ in $G$, there exists an expression
\begin{equation}
W=_F\prod\limits_{i=1}^k f_i^{-1}R_i^{\pm 1}f_i \label{prod}
\end{equation}
with the equality in the group $F$, where $R_i\in \mathcal R$ and
$f_i\in F $ for $i=1, \ldots , k$. The smallest possible number
$k$ in a representation of the form (\ref{prod}) is called the
{\it relative area} of $W$ and is denoted by $Area^{rel}(W)$.

\begin{defn}\label{IP}
We say that a function $f:\mathbb N\to \mathbb N$ is a {\it
relative isoperimetric function} of $G$ with respect to $\Hl $
(associated to the relative presentation (\ref{G})) if for any
$n\in \mathbb N$ and any word $W$ in $X\cup \mathcal H$ of length
at most $n$ representing the identity in the group $G$, we have
$Area^{rel} (W)\le f(n).$ The smallest relative isoperimetric
function of (\ref{G}) is called the {\it relative Dehn function}
of $G$ with respect to $\Hl $ and is denoted by
$\delta^{rel}_{G,\, \Hl }$ (or simply by $\delta ^{rel}$ when the
group $G$ and the collection of subgroups are fixed).
\end{defn}

Observe that $\delta ^{rel} (n)$ is not always well--defined,
since the number of words of bounded relative length can be
infinite. Indeed consider the group
$$G=\langle a,b \; |\; [a,b]=1 \rangle\cong \mathbb Z\times
\mathbb Z $$ and the cyclic subgroup $H$ generated by $a$. Clearly
$X=\{ b^{\pm 1}\} $ is a relative generating set of $G$ with
respect to $H$. It is easy to see that the word $W_n=[a^n,b]$ has
length $4$ as a word over $X\cup (H\setminus \{ 1\} )$ for every
$n$ since $a^n$ can be regarded as a letter from the alphabet
$H\setminus \{ 1\} $. On the other hand $Area^{rel} (W_n)$ growths
linearly as $n\to\infty $. Thus we do not have any bound on
$Area^{rel} (W_n)$ in terms of $\| W_n\| $ in this case.

However if $G$ is finitely presented relative to $\Hl $ and
$\delta ^{rel}$ is well--defined, it is independent of the choice
of the finite relative presentation up to the following
equivalence relation \cite[Theorem 2.32]{RHG}. For two functions
$f,g:\mathbb N\to \mathbb N$, we write $f\preceq g$ if there are
positive constants $A,B,C$ such that $f(n)\le Ag(Bn)+Cn.$ One says
that $f$ and $g$ are {\it equivalent} if $f\preceq g$ and
$g\preceq f$.

Recall that a function ${f : {\mathbb N}\to {\mathbb N}}$ is said
to be {\it superadditive} if
$$f(a+b)\ge f(a)+f(b)$$ for any $a,b\in {\mathbb N}$. Given an
arbitrary function  ${f : {\mathbb N}\to {\mathbb N}}$, the {\it
superadditive closure} of $f$ is defined to be
\begin{equation}
\overline f(n)=\max_{i=1,\dots , n}\left(\max_{a_1+\dots + a_i=n,\
a_i\in {\mathbb N}} \left(f(a_1)+\dots +f(a_i)\right)\right)
\end{equation}
In fact, $\overline f$ is the smallest superadditive function such
that $\overline f(n)\ge f(n)$ for all $n$. The main results of our
paper are the following two theorems.

\begin{thm}\label{HNN}
Suppose that a group $H$ is finitely presented with respect to a
collection of subgroups $\Hl \cup \{ K\} $ and the corresponding
relative Dehn function $\gamma $ is well--defined. Assume also
that $K$ is finitely generated and for some $\nu \in \Lambda $,
there exists a monomorphism $\iota :K\to H_{\nu }$. Then the
HNN--extension
\begin{equation}\label{HNN-pres}
 G=\langle H,t\; |\; t^{-1}kt=\iota (k),\; k\in K\rangle
\end{equation}
is finitely presented relative to $\Hl $. Moreover, the relative
Dehn function $\delta $ of $G$ with respect to $\Hl $ is
well--defined and satisfies the inequality
$$\delta \preceq \overline \gamma \circ \overline \gamma.$$
\end{thm}

Similarly for amalgamated products, we have

\begin{thm}\label{Am}
Suppose that a group $A$ (respectively $B$) is finitely presented
relative to a collection of subgroups $\Am \cup \{ K\} $
(respectively $\Bn $) and the corresponding relative Dehn function
$\gamma _1$ (respectively $\gamma _2$) is well--defined. Assume in
addition that $K$ is finitely generated and for some $\eta\in N$,
there is a monomorphism $\xi : K\to B_\eta $. Then the amalgamated
product $A\ast _{K=\xi (K)} B$ is finitely presented relative to
the collection $\Am \cup\Bn $. Moreover, the relative Dehn
function $\delta $ of $A\ast _{K=\xi (K)} B$ with respect to the
collection $\Am \cup\Bn$ is well--defined and satisfies the
inequality
$$\delta \preceq \overline \gamma \circ \overline \gamma,$$ where $\gamma $ is defined by
\begin{equation} \label{g}
\gamma (n)=\max \{ \gamma _1(n),\, \gamma _2(n)\}
\end{equation}
for all $n\in \mathbb N$.
\end{thm}

In fact Theorem \ref{Am} can be derived from Theorem \ref{HNN} via
an easy observation concerning the behavior of relative Dehn
functions under taking retracts (see Section 3).

Recall that a group $G$ is {\it hyperbolic relative to a
collection of subgroups} $\Hl $ if the relative Dehn function of
$G$ with respect to $\Hl $ is linear \cite{RHG}. Thus we obtain
the following corollaries of Theorem \ref{HNN} and Theorem
\ref{Am}. For finitely generated groups, these results were
obtained by Dahmani in \cite{Dah} and used to prove that limits
groups introduced by Sela \cite{Sela} are hyperbolic relative to
certain collections of maximal abelian subgroups.

\begin{cor}\label{HNN-hyp}
In the notation of Theorem \ref{HNN}, if $H$ is hyperbolic
relative to $\Hl \cup \{ K\} $, then $G$ is hyperbolic relative to
$\Hl $.
\end{cor}

\begin{cor}
In the notation of Theorem \ref{Am}, if $A$ and $B$ are hyperbolic
relative to $\Am \cup \{ K\} $ and $\Bn$ respectively, then $A\ast
_{K=\xi (K)} B$ is hyperbolic relative to $\Am \cup\Bn $.
\end{cor}

\begin{ex}
Recall that for any quasi--convex subgroups $Q_1, Q_2$ of an
ordinary hyperbolic group $H$, $H$ is hyperbolic relative to $\{
Q_1, Q_2\} $ whenever $|Q_i^h\cap Q_i|<\infty $ for any $i=1,2$
and $h\notin Q_i$, and $|Q_1^g\cap Q_2|< \infty $ for any $g\in H$
\cite{Bow,ESBG}. Suppose that there is a monomorphism $\iota
\colon Q_1\to Q_2$. Then the corresponding HNN--extension $H\ast
_\iota $ is hyperbolic relative to $Q_2$ by Corollary
\ref{HNN-hyp}. As $Q_2$ is hyperbolic itself (see \cite{Gro}),
$H\ast _\iota $ is a hyperbolic group \cite{F}.

Similarly suppose that $Q_1$ and $Q_2$ are quasi--convex subgroups
of hyperbolic groups $H_1$ and $H_2$ respectively, $|Q_1^g\cap
Q_1|<\infty $ for any $g\notin Q_1$, and $|Q_2^h\cap Q_2|<\infty $
for any $h\notin Q_2$. Assume that there is a monomorphism $\xi
\colon Q_1\to Q_2$. Then the amalgamated product $H_1\ast
_{Q_1=\xi (Q_1)}H_2 $ is hyperbolic relative to $Q_2$. In
particular, $H_1\ast _{Q_1=\xi (Q_1)}H_2 $ is hyperbolic. These
results are well known and were proved independently by many
authors (see \cite{BF,KM,MO}).
\end{ex}


\section{Preliminaries}


\noindent {\bf Word metrics and Cayley graphs. } Let $G$ be a
group generated by a symmetrized set $\mathcal A$. Recall that the
{\it Cayley graph} $\Ga $ of a group $G$ with respect to the set
of generators $\mathcal A$ is an oriented labelled 1--complex with
the vertex set $V(\Ga )=G$ and the edge set $E(\Ga )=G\times
\mathcal A$. An edge $e=(g,a)$ goes from the vertex $g$ to the
vertex $ga$ and has the label $\phi (e)=a$. As usual, we denote
the origin and the terminus of the edge $e$, i.e., the vertices
$g$ and $ga$, by $e_-$ and $e_+$ respectively. Given a
combinatorial path $p=e_1e_2\ldots e_k$ in the Cayley graph $\Ga
$, where $e_1, e_2, \ldots , e_k\in E(\Ga )$, we denote by $\phi
(p)$ its label. By definition, $\phi (p)=\phi (e_1)\phi(e_2)\ldots
\phi (e_k).$ We also denote by $p_-=(e_1)_-$ and $p_+=(e_k)_+$ the
origin and the terminus of $p$ respectively. The length $l(p)$ of
$p$ is, by definition, the number of edges of $p$.

Associated to $\mathcal A$ is the so--called {\it word metric} on
$G$. More precisely, the length $|g|_\mathcal A$ of an element
$g\in G$ is defined to be the length of a shortest word in
$\mathcal A^{\pm 1}$ representing $g$ in $G$. This defines a
metric on $G$ by $dist_\mathcal A(f,g)=|f^{-1}g|_\mathcal A$. We
also denote by $dist _\mathcal A$ the natural extension of the
word metric to the Cayley graph $\Ga $.

\vspace{2mm}

\noindent {\bf Van Kampen Diagrams. } Recall that a {\it van
Kampen diagram} $\Delta $ over a presentation
\begin{equation}
G=\langle \mathcal A\; | \; \mathcal O\rangle \label{ZP}
\end{equation}
is a finite oriented connected simply--connected 2--complex
endowed with a labelling function $\phi : E(\Delta )\to \mathcal
A$, where $E(\Delta ) $ denotes the set of oriented edges of
$\Delta $, such that $\phi (e^{-1})=(\phi (e))^{-1}$. Labels and
lengths of paths are defined as in the case of Cayley graphs.
Given a cell $\Pi $ of $\Delta $, we denote by $\partial \Pi$ the
boundary of $\Pi $; similarly, $\partial \Delta $ denotes the
boundary of $\Delta $. The labels of $\partial \Pi $ and $\partial
\Delta $ are defined up to a cyclic permutation. An additional
requirement is that for any cell $\Pi $ of $\Delta $, the boundary
label $\phi (\partial \Pi)$ is equal to a cyclic permutation of a
word $P^{\pm 1}$, where $P\in \mathcal O$. Sometimes it is
convenient to use the notion of $0$--refinement in order to assume
diagrams to be homeomorphic to a disc. We do not explain here this
notion and refer the interested reader to \cite[Ch. 4]{Ols-book}.

The van Kampen lemma states that a word $W$ over the alphabet
$\mathcal A$ represents the identity in the group given by
(\ref{ZP}) if and only if there exists a simply--connected planar
diagram $\Delta $ over (\ref{ZP}) such that with boundary label
$W$ \cite[Ch. 5, Theorem 1.1]{LS}.

\vspace{2mm}

\noindent {\bf $H_\lambda $--components. } For a group $G$ and a
collection of subgroups $\Hl $ of $G$, we denote by $\GG $ the
Cayley graph of $G$ with respect to the generating set $\mathcal
A=X\cup \mathcal H$, where $\mathcal H$ is defined by (\ref{H}).
We also fix a relative presentation (\ref{G}) of $G$ with respect
to $\Hl $.

Let $q$ be a path in $\GG $. A subpath $p$ of $q$ is called an
{\it $H_\lambda $--component} for some $\lambda \in \Lambda $ (or
simply a {\it component}) of $q$, if the label of $p$ is a word in
the alphabet $H_\lambda\setminus \{ 1\} $ and $p$ is not contained
in a bigger subpath of $q$ with this property. Two $H_\lambda
$--components $p_1, p_2$ of a path $q$ in $\GG $ are called {\it
connected} if there exists a path $c$ in $\GG $ that connects some
vertex of $p_1$ to some vertex of $p_2$ and ${\phi (c)}$ is a word
consisting of letters from $ H_\lambda\setminus\{ 1\} $. In
algebraic terms this means that all vertices of $p_1$ and $p_2$
belong to the same coset $gH_\lambda $ for a certain $g\in G$.
Note that we can always assume $c$ to have length at most $1$, as
every nontrivial element of $H_\lambda $ is included in the set of
generators.  An $H_\lambda $--component $p$ of a path $q$ is
called {\it isolated } if no distinct $H_\lambda $--component of
$q$ is connected to $p$.

The next lemma is a simplification of Lemma 2.27 from \cite{RHG}.

\begin{lem}\label{Omega}
Suppose that the relative presentation (\ref{G}) is finite and the
corresponding relative Dehn function $\delta $ is well--defined.
Then for each subgroup $H_\lambda $ there exists a constant $C>0$
and a finite subset $\Omega _{H_\lambda } \subseteq H_\lambda $,
such that the following condition holds. Let $q$ be a cycle in $\G
$, $p_1, \ldots , p_k$ a set of isolated $H_\lambda $--components
of $q$, $g_1, \ldots , g_k$ the elements of $G$ represented by the
labels of $p_1, \ldots , p_k$ respectively. Then for any $i=1,
\ldots , k$, $g_i$ belongs to the subgroup $\langle \Omega
_{H_\lambda } \rangle \le G$ and the lengths of $g_i$ with respect
to $\Omega _{H_\lambda } $ satisfy the inequality
$$ \sum\limits_{i=1}^k |g_i|_{\Omega _{H_\lambda } }\le C\delta (l(q)).$$
\end{lem}

\vspace{2mm}

\noindent {\bf Conventions and notation.} Given a word $W$ in an
alphabet $\mathcal A$, we denote by $\| W\| $ its length. We also
write $W\equiv V$ to express letter--for--letter equality of words
$W$ and $V$. Finally for elements $g$, $t$ of a group $G$, $g^t$
denotes the element $t^{-1}gt$. Similarly $H^t$ denotes $t^{-1}Ht$
for a subgroup $H$ of $G$.


\section{Diagrams over HNN--extensions}


To deal with van Kampen diagrams over HNN--extensions, we need
certain standard technical tools. The first one is a
\emph{$t$-band} (other people call them corridors or strips). This
notion goes back to the paper \cite{MS}; here we describe it
shortly. Suppose that one has a finite presentation of the form
\begin{equation}\label{HNN-type}
\langle \mathcal A,\, t \, \mid \, \mathcal{T}\rangle ,
\end{equation}
where $t\notin \mathcal{A}$ and the only relators involving $t$
are of the form
\begin{equation}\label{rel}
t^{-1}x_ity_i
\end{equation}
where $x_i, y_i$ are some words over $\mathcal{A}$. If the
boundary label $\partial \pi$ of a cell $\pi $ of a van Kampen
diagram over (\ref{HNN-type}) is of the form~(\ref{rel}), we call
$\pi $ a $t$-\emph{cell}. The subpath of $\partial \pi $ having
the label $x_i$ (respectively $y_i$) is called a \emph{bottom}
(respectively \emph{top}) of $\pi $. A $t$-band is a sequence of
pairwise distinct $t$-cells $\tau =(\pi _1, \ldots , \pi _n)$ in a
van Kampen diagram over (\ref{HNN-type}) such that each two
consecutive cells in this sequence have a common edge labelled
$t$. The path formed by tops (respectively bottoms) of $\pi _1,
\ldots , \pi _n$ is called a {\it top} (respectively {\it bottom})
of the $t$-band $\tau $ and is denoted by $\t (\tau )$
(respectively $\b (\tau )$). The \emph{length} of a $t$-band is
the number of its 2-cells. A $t$-band $\pi _1, \ldots , \pi _n$ is
\emph{maximal} if it is not contained in another $t$-band of
greater length. In what follows we always assume all $t$--bands
under consideration to be maximal. Finally a (maximal) $t$-band
$\tau =(\pi _1, \ldots , \pi _n)$ is called a $t$-annulus if it
has no common edges with $\partial \Delta $ labelled $t$.

We will use the following obvious information about $t$-bands.

$\bullet $ Distinct $t$-bands have no common 2-cells.

$\bullet $ The labels of the top and the bottom of any $t$-band
contain no $t^{\pm 1}$.

\begin{figure}

\begin{picture}(105,52)(18,56)
\put(74.5,78.38){\oval(61,38.75)[]}
\qbezier(44,88.89)(57.68,80.27)(44,69.27)
\qbezier(63.18,97.51)(73.14,74.03)(53.96,58.87)

\qbezier(88,97.66)(89.04,84.14)(104.95,86.07)

\qbezier(74.03,58.87)(88.6,74.62)(104.95,77.89)

\qbezier(83.99,58.87)(91.27,72.76)(97.96,59.01)

\put(81,82){\makebox(0,0)[cc]{$D $}}
\put(68,76.56){\makebox(0,0)[cc]{$\tau $}}

\put(88.5,73){\makebox(0,0)[cc]{$\sigma $}}

\put(71.2,86.03){\vector(-1,-1){.07}}
\multiput(86,102.87)(-.0168555905,-.0192635321){860}{\line(0,-1){.0192635321}}
\put(63.62,70.76){\vector(4,1){.07}}
\multiput(37.01,64.07)(.067023667,.016849525){397}{\line(1,0){.067023667}}
\put(32,61){\makebox(0,0)[cc]{$\t (\tau )$ }}
\put(85.26,105){\makebox(0,0)[cc]{$\b (\tau )$}}

\put(56.3,58.97){\vector(-1,0){.07}} \put(56.3,55){$t$}

\put(64.5,97.75){\vector(-1,0){.07}} \put(64.5,99){$t$}

\put(41,102){\vector(1,-1){13}} \put(37,103){$t$--annulus}

\put(58,85){\circle{10}} \put(75,72){\circle{7}}
\put(90,82){\circle{4}}

\thicklines \qbezier(44,84.88)(51.43,79.9)(44,74.03)
\qbezier(67.04,97.66)(77,75.22)(59.01,58.72)

\put(75,72){\circle{4}} \put(90,82){\circle{7}}
\put(58,85){\circle{6}}

\qbezier(90.97,97.51)(92.31,89.26)(104.95,89.04)
\qbezier(70.01,59.01)(83.76,77.3)(104.95,81.91)
\qbezier(86.96,59.01)(91.2,66.82)(95.14,58.87)

\linethickness{0.1pt}
\multiput(86.07,97.66)(-.0168555905,-.0203872381){926}{\line(0,-1){.0203872381}}
\multiput(88.3,94.99)(-.016851041,-.019860156){494}{\line(0,-1){.019860156}}
\multiput(89.64,91.57)(-.016844222,-.021055277){353}{\line(0,-1){.021055277}}
\multiput(91.57,88.74)(-.0168605459,-.0213424631){1393}{\line(0,-1){.0213424631}}
\multiput(79.83,79.23)(-.0168523285,-.0214628711){935}{\line(0,-1){.0214628711}}
\multiput(76.7,80.72)(-.016860142,-.020232171){529}{\line(0,-1){.020232171}}
\multiput(94.54,86.96)(-.0168681141,-.0212959941){705}{\line(0,-1){.0212959941}}
\multiput(97.81,86.07)(-.016868114,-.021436562){423}{\line(0,-1){.021436562}}
\multiput(101.08,85.62)(-.016818095,-.022243286){274}{\line(0,-1){.022243286}}
\multiput(104.2,85.92)(-.01684215,-.02416482){203}{\line(0,-1){.02416482}}
\multiput(76.11,97.51)(-.016865332,-.019610852){379}{\line(0,-1){.019610852}}
\multiput(71.95,97.66)(-.01685725,-.01915596){194}{\line(0,-1){.01915596}}
\multiput(62.28,61.84)(.02477504,-.01681164){168}{\line(1,0){.02477504}}
\multiput(64.66,64.66)(.023807268,-.016839287){256}{\line(1,0){.023807268}}
\multiput(66.89,68.08)(.023870844,-.016818095){274}{\line(1,0){.023870844}}
\multiput(68.53,71.65)(.022186605,-.01686182){335}{\line(1,0){.022186605}}
\multiput(69.87,75.51)(.021601884,-.016842147){406}{\line(1,0){.021601884}}
\multiput(70.61,80.57)(.020586456,-.016843464){556}{\line(1,0){.020586456}}
\multiput(71.5,84.73)(.02046471,-.016866519){661}{\line(1,0){.02046471}}
\multiput(74.33,88.15)(.02123575,-.01681164){168}{\line(1,0){.02123575}}
\multiput(83.39,81.01)(.020558014,-.016868114){282}{\line(1,0){.020558014}}
\multiput(77.6,91.42)(.0208799973,-.0168572455){776}{\line(1,0){.0208799973}}
\multiput(68.68,93.8)(.01748827,-.01681564){221}{\line(1,0){.01748827}}
\multiput(70.76,97.66)(.019258333,-.016851041){247}{\line(1,0){.019258333}}
\multiput(80.87,94.69)(.0219349768,-.0168591144){820}{\line(1,0){.0219349768}}
\multiput(84.14,97.51)(.02222091,-.01685725){194}{\line(1,0){.02222091}}
\multiput(99.89,85.77)(.02637343,-.01678309){186}{\line(1,0){.02637343}}
\end{picture}
\caption{Structure of diagrams over HNN--extensions; bottoms of
$t$--bands are marked in bold.} \label{pic0}
\end{figure}
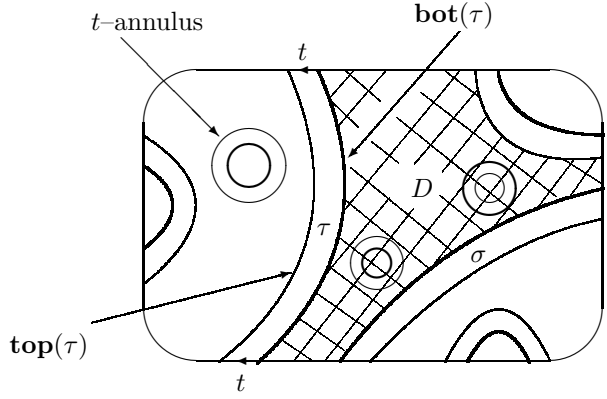

Furthermore, $\partial \Delta $ and non--annular  $t$-bands bound
{\it domains} of $\Delta $, which are maximal connected disjoint
subdiagrams of $\Delta$ containing no $t$--cells except for those
from $t$--annuli. The set of all domains of $\Delta $ is denoted
by $\mathcal D(\Delta )$. By $T(D)$ we denote the set of all
non--annular $t$--bands $\tau $ such that $\b (\tau )$ lies on the
boundary of a domain $D$. For example, $T(D)=\{ \tau , \sigma\} $
in Fig. \ref{pic0}.

In the rest of the paper we use the notation of Theorem \ref{HNN}.
Let us fix a finite presentation
\begin{equation}\label{1}
H=\langle X, K, H_\lambda , \lambda\in \Lambda \; |\; \mathcal R
\rangle
\end{equation}
of $H$ with respect to $\Hl \cup \{ K\} $.  We denote by $Y$ an
arbitrary symmetrized finite generating set of $K$. For a word
$R\in \mathcal R$, let us denote by $R^\prime $ the word obtained
from $R$ by replacing each letter $k\in K$ with a shortest word in
$Y$ representing the element $k$ in $H$. Set $\mathcal R^\prime
=\{ R^\prime , \, R\in \mathcal R\} $, $\mathcal T=\{
t^{-1}yt=\iota (y),\, y\in Y \} $ (here $\iota (y)\in H_\nu $ is
regarded as a letter), and $\widetilde{\mathcal R}=\mathcal
R^\prime \cup \mathcal T$. Obviously the HNN--extension
(\ref{HNN-pres}) has the finite relative presentation
\begin{equation}\label{2}
G=\langle X\cup Y\cup\{ t\} , H_\lambda , \lambda\in \Lambda \;
|\; \widetilde{ \mathcal R} \rangle .
\end{equation}

Associated to the relative presentations (\ref{1}) and (\ref{2})
are ordinary (non--relative) presentations

\begin{equation}\label{p1}
H=\langle X\cup\mathcal H\cup ( K\setminus\{ 1\} ) \; |\; \mathcal
R\cup \mathcal S\cup\mathcal L \rangle
\end{equation}
and
\begin{equation}\label{p2} G=\langle X\cup Y\cup\{t \}\cup\mathcal H
\; |\; \widetilde{ \mathcal R}\cup \mathcal S \rangle ,
\end{equation}
where $\mathcal H$ is defined by (\ref{H}), $\mathcal
S=\bigcup\limits_{\lambda \in \Lambda } S_\lambda $ and $S_\lambda
$ (respectively $\mathcal L$) is the set of all words in the
alphabet $H_\lambda \setminus \{ 1\} $ (respectively $K\setminus
\{ 1\} $) representing $1$ in $H$.

\begin{defn}
Consider two $t$--bands $\tau _1, \tau _2\in T(D)$, where $D$ is a
domain of a diagram $\Delta $ over (\ref{p2}). The boundary of $D$
can be decomposed as $$\partial D=\b (\tau _1) a\b(\tau _2)b$$ for
certain subpaths $a,b$ of $\partial D$. We say that $\tau_1 $ and
$\tau _2$ are {\it $K$--connected} if $\phi (a)$ (or,
equivalently, $\phi (b)$) represents an element of $K$ in $G$.
\end{defn}

The following lemma plays the key role in the proof of Theorem
\ref{HNN}.

\begin{lem}\label{mainl}
Suppose that $W$ is a word in the alphabet $X\cup Y\cup\{ t\} \cup
\mathcal H$ such that $W=1$ in $G$.  Then there is a van Kampen
diagram $\Delta $ over (\ref{p2}) with boundary label $W$ such
that the following conditions hold.
\begin{enumerate}
\item $\Delta $ contains no $K$--connected $t$--bands. \item Any
non--annular $t$--band $\tau $ in $\Delta $ has length $|k(\tau
)|_Y$, where $k(\tau )$ denotes the element of $K$ represented by
$\phi (\b (\tau ))$.
\end{enumerate}
\end{lem}

\begin{figure}
\unitlength 0.94mm 
\linethickness{0.4pt}
\ifx\plotpoint\undefined\newsavebox{\plotpoint}\fi 
\begin{picture}(138.99,75)(8,10)

\put(15.25,62.38){\circle*{1}} \put(77.09,61.93){\circle*{1}}
\put(53.88,62.25){\circle*{1}} \put(135.93,62.55){\circle*{1}}
\put(15.25,25){\circle*{1}} \put(77.09,24.55){\circle*{1}}
\put(53.88,24.75){\circle*{1}} \put(135.93,25.05){\circle*{1}}
\put(53.88,29.88){\circle*{1}} \put(135.93,30.17){\circle*{1}}
\put(54,56.63){\circle*{1}} \put(136.06,56.92){\circle*{1}}
\put(15.38,56.75){\circle*{1}} \put(77.21,56.3){\circle*{1}}
\put(15.13,29.88){\circle*{1}} \put(76.96,29.43){\circle*{1}}

\put(86.96,24.97){\circle*{1}} \put(86.96,62.43){\circle*{1}}
\put(125.76,62.73){\circle*{1}} \put(125.76,24.97){\circle*{1}}

\qbezier(53.63,17.88)(48.19,17.63)(48,11.88)
\qbezier(125.42,18.02)(119.98,17.77)(119.8,12.02)
\qbezier(53.88,19.88)(45.81,19.13)(46,11.38)
\qbezier(125.67,20.02)(117.61,19.27)(117.8,11.52)
\qbezier(23.97,11.35)(27.96,16.71)(31.95,11.56)
\qbezier(95.76,11.5)(99.76,16.86)(103.75,11.71)
\qbezier(21.97,11.46)(28.8,20.71)(33.95,11.46)
\qbezier(93.77,11.61)(100.6,20.86)(105.75,11.61)
\qbezier(15.24,20.92)(34.74,24.81)(35.95,11.46)
\qbezier(87.04,21.07)(106.54,24.96)(107.75,11.61)
\qbezier(15.24,22.91)(36.84,26.91)(38.05,11.56)
\qbezier(87.04,23.06)(108.64,27.06)(109.85,11.71)
\qbezier(27.85,75.79)(34.06,63.38)(41.1,75.79)
\qbezier(99.65,75.93)(105.85,63.53)(112.9,75.93)
\qbezier(25.96,75.68)(34.48,59.7)(42.99,75.68)
\qbezier(97.76,75.83)(106.27,59.85)(114.79,75.83)

\qbezier(15.24,53.92)(31.9,39.89)(15.14,32.79)
\qbezier(77.08,53.48)(93.74,39.44)(76.97,32.35)
\qbezier(15.03,52.03)(28.22,39.78)(15.14,34.9)
\qbezier(76.87,51.58)(90.06,39.34)(76.97,34.45)
\qbezier(15.14,36.05)(19.34,38.05)(15.14,40.89)
\qbezier(76.97,35.61)(81.18,37.6)(76.97,40.44)
\qbezier(15.14,37.95)(16.13,38.47)(15.24,39)
\qbezier(76.97,37.5)(77.97,38.03)(77.08,38.55)
\qbezier(15.24,49.93)(21.13,45.99)(15.24,43.94)
\qbezier(77.08,49.48)(82.97,45.54)(77.08,43.49)
\qbezier(15.14,47.93)(17.24,46.99)(15.14,46.04)
\qbezier(76.97,47.49)(79.08,46.54)(76.97,45.59)
\qbezier(54.03,52.98)(42.52,40.94)(53.92,39.21)
\qbezier(136.08,53.27)(124.57,41.24)(135.98,39.5)
\qbezier(54.13,50.98)(46.51,41.41)(53.82,41.1)
\qbezier(136.19,51.28)(128.57,41.71)(135.87,41.4)

\put(34.34,78.7){\makebox(0,0)[cc]{$U_2$}}
\put(41.62,7.88){\makebox(0,0)[cc]{$U_1$}}
\put(11.45,42.37){\makebox(0,0)[cc]{$V_1$}}
\put(12.49,27.5){\makebox(0,0)[cc]{$t$}}
\put(12.49,59){\makebox(0,0)[cc]{$t$}}
\put(56,28){\makebox(0,0)[cc]{$t$}}

\put(55.89,59.76){\makebox(0,0)[cc]{$t$}}
\put(34.19,55.15){\makebox(0,0)[cc]{$\tau _1$}}
\put(30.03,32.85){\makebox(0,0)[cc]{$\tau _2$}}
\put(74.47,58.87){\makebox(0,0)[cc]{$t$}}
\put(74.03,27.05){\makebox(0,0)[cc]{$t$}}
\put(89.5,44.15){\makebox(0,0)[cc]{$B_1$}}
\put(80,21){\makebox(0,0)[cc]{$\Delta _1$}}
\put(106.14,30.47){\makebox(0,0)[cc]{$\Delta _0$}}
\put(34.49,44){\makebox(0,0)[cc]{$\Delta $}}
\put(134,21){\makebox(0,0)[cc]{$\Delta _2$}}
\put(124,43.41){\makebox(0,0)[cc]{$B_2$}}

\put(71.06,41.92){\vector(1,0){.07}}\put(59.01,41.92){\line(1,0){12.041}}

\qbezier(15.13,62.5)(36,53)(53.88,62.5)
\qbezier(15.13,56.75)(35.94,48.63)(54,56.5)
\qbezier(15,29.75)(24.06,41)(53.88,29.75)
\qbezier(15.13,25.13)(25.31,36.75)(53.75,24.88)

\put(92.7,44){\vector(0,1){.07}}\qbezier(77,29.43)(108.22,45.12)(77,56.34)
\put(120.5,43.26){\vector(0,-1){.07}}\qbezier(135.87,56.78)(104.95,43.03)(135.87,30.18)


\put(15.14,43.25){\vector(0,1){.07}}\put(15.14,29.85){\line(0,1){26.804}}
\put(76.97,42.81){\vector(0,1){.07}}\put(76.97,29.41){\line(0,1){26.804}}
\put(53.92,43.31){\vector(0,-1){.07}}\multiput(54.03,56.55)(-.03003,-3.78403){7}{\line(0,-1){3.78403}}
\put(135.98,43.6){\vector(0,-1){.07}}\multiput(136.08,56.85)(-.03003,-3.78403){7}{\line(0,-1){3.78403}}

\put(102.97,43.85){\vector(0,1){.07}}\qbezier(86.96,24.97)(118.77,44)(86.96,62.43)

\put(110.25,42.81){\vector(0,-1){.07}}\qbezier(125.76,62.73)(94.54,41.77)(125.76,24.97)
\put(97.61,43.85){\vector(0,1){.07}}\qbezier(77,24.97)(117.88,44.3)(77.3,61.84)

\put(115.9,43.11){\vector(0,-1){.07}}\qbezier(136.16,62.73)(95.43,42.44)(136.16,24.82)

\put(35.08,75.85){\vector(1,0){.07}}
\put(41.7,11.52){\vector(-1,0){.07}}
\put(105.91,75.97){\vector(1,0){.07}}
\put(114.39,11.55){\vector(-1,0){.07}}

\put(34.56,62.4){\oval(38.71,26.87)[t]}
\put(106.36,62.55){\oval(38.71,26.87)[t]}
\put(34.56,24.93){\oval(38.71,26.87)[b]}
\put(106.36,25.07){\oval(38.71,26.87)[b]}

\put(15.16,56.78){\line(0,1){5.649}}
\put(15.01,29.88){\line(0,-1){4.608}}
\put(53.81,29.73){\line(0,-1){4.757}}
\put(53.81,62.14){\line(0,-1){5.351}}
\put(77.15,61.84){\line(0,-1){5.351}}
\put(77,29.43){\line(0,-1){4.608}}
\put(135.87,62.43){\line(0,-1){5.649}}
\put(135.87,30.03){\line(0,-1){4.905}}

\put(15.01,26.61){\vector(0,-1){.07}}
\put(53.81,26.76){\vector(0,-1){.07}}
\put(53.81,59.76){\vector(0,1){.07}}
\put(15.16,59.91){\vector(0,1){.07}}
\put(77.15,59.61){\vector(0,1){.07}}
\put(77,26.31){\vector(0,-1){.07}}
\put(135.87,60.5){\vector(0,1){.07}}
\put(135.87,26.76){\vector(0,-1){.07}}
\put(138.8,59.76){\makebox(0,0)[cc]{$t$}}
\put(138.8,27.8){\makebox(0,0)[cc]{$t$}}
\put(139.2,45.93){\makebox(0,0)[cc]{$V_2$}}
\put(106.58,78.7){\makebox(0,0)[cc]{$U_2$}}
\put(114.91,7.88){\makebox(0,0)[cc]{$U_1$}}
\put(100.5,45){\makebox(0,0)[cc]{$T_1$}}
\put(113.53,42){\makebox(0,0)[cc]{$T_2$}} \thinlines

\put(57,45.34){\makebox(0,0)[cc]{$V_2$}}
\put(73.88,39){\makebox(0,0)[cc]{$V_1$}}
\end{picture}

  \caption{Rebuilding of diagrams containing $K$--connected $t$--bands.}\label{fig1}
\end{figure}
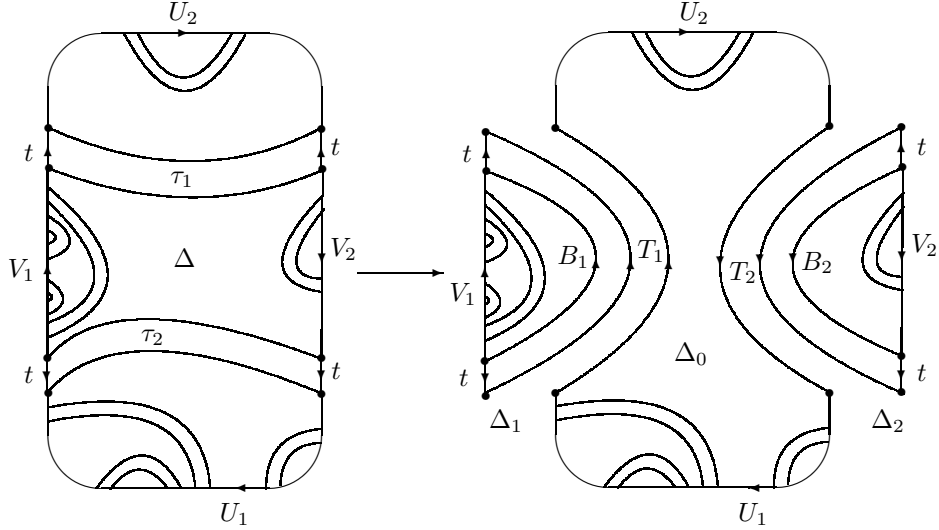

\begin{proof}
To prove the first assertion of the lemma we proceed by induction
on the number $n$ of appearances of the letters $t^{\pm 1} $ in
$W$. If $n=0$, the lemma is obvious, so we assume $n>0$. Let
$\Delta $ be an arbitrary diagram over (\ref{p2}) with boundary
label $W$. Suppose that for some $D\in \mathcal D(\Delta )$, there
are two $K$--connected $t$--bands $\tau _1, \tau _2\in T(D)$. Then
the boundary of $\Delta $ is decomposed as
$$\partial \Delta =u_1e_1^{-1}v_1e_2u_2e_3^{-1}v_2e_4,$$ where
$e_2, e_3$ (respectively $e_1, e_4$) are common edges of $\partial
\Delta $ and $\tau _1$ (respectively $\tau _2$) labelled $t$. Let
$U_i$, $V_i$ denote the labels of $u_i$ and $v_i$, $i=1,2$ (see
Fig.\ref{fig1}). We are going to rebuild the diagram $\Delta $ as
follows.

Since $\tau _1$ and $\tau _2$ are $K$--connected, $V_i$ represents
an element of $K$ in $G$ for $i=1,2$. Let $B_i$ be a word in $Y$
such that $B_i=V_i$ in $G$. Further let $T_i$ be the word in the
alphabet $H_\nu \setminus \{ 1\} $ obtained from $B_i$ by
replacing each letter $y\in Y$ with the letter from $H_\nu
\setminus \{ 1\} $ representing the element $y^t$ in $G$. For
$i=1,2$, we denote by $\Delta_i$ the van Kampen diagram over
(\ref{p2}) with boundary label
$$\phi (\Delta _i)\equiv t^{-1}V_itT_i^{-1}$$ that is obtained by
gluing the $t$--band with the bottom labelled $B_i$ to a diagram
over (\ref{p2}) corresponding to the equality $V_i=B_i$ in $G$.
Clearly $$U_1T_1U_2T_2=U_1t^{-1}V_1tU_2t^{-1}V_2t=1$$ in $G$.
Denote by $\Delta _0$ a diagram over (\ref{p2}) with boundary
label
$$\phi (\Delta _0)\equiv U_1T_1U_2T_2.$$

We now glue the diagrams $\Delta _0$, $\Delta _1$, $\Delta _2$ in
the obvious way (see Fig. \ref{fig1}) and denote the obtained
diagram by $\Xi $. It is clear that any domain $D$ of $\Xi $ is a
domain in $\Delta _j$ for a certain $j\in \{ 0,1,2\} $ and two
$t$--bands $\tau _1, \tau _2\in T(D)$ are $K$--connected in $\Xi $
if and only if they are $K$--connected in $\Delta _j$. Since the
number of appearances of $t^{\pm 1}$ in the boundary labels of
each of the diagrams $\Delta _0$, $\Delta _1$, $\Delta _2$ is at
most $(n-2)$, we may assume that $\Delta _0$, $\Delta _1$, and
$\Delta _2$ contain no $K$--connected $t$--bands. Therefore no
$t$--bands of $\Xi $ are $K$--connected.

\begin{figure}
  \unitlength 1mm 
\linethickness{0.4pt}
\ifx\plotpoint\undefined\newsavebox{\plotpoint}\fi 
\begin{picture}(77.38,68)(3,18)
\put(20.59,53.21){\oval(24.93,27.93)[t]}
\put(64.47,70.21){\oval(24.93,27.93)[t]}
\put(20.59,47.55){\oval(24.93,27.93)[b]}
\put(64.84,30.93){\oval(24.93,27.93)[b]}
\put(20.5,53.31){\vector(1,0){.07}}\multiput(8,53.38)(6.25,-.0313){4}{\line(1,0){6.25}}
\put(64.75,53.31){\vector(1,0){.07}}\multiput(52.25,53.38)(6.25,-.0313){4}{\line(1,0){6.25}}
\put(20.56,47.44){\vector(1,0){.07}}\multiput(8.13,47.5)(6.2188,-.0313){4}{\line(1,0){6.2188}}
\put(64.81,47.44){\vector(1,0){.07}}\multiput(52.38,47.5)(6.2188,-.0313){4}{\line(1,0){6.2188}}
\put(8.13,47.38){\line(0,1){6.125}}
\put(52.38,47.38){\line(0,1){6.125}}
\put(33,47.38){\line(0,1){5.875}}
\put(77.25,47.38){\line(0,1){5.875}}
\put(8.1,50.44){\vector(0,1){.07}}
\put(52.35,50.44){\vector(0,1){.07}}
\put(33,50.44){\vector(0,1){.07}}
\put(77.25,50.44){\vector(0,1){.07}}
\put(64.44,70.25){\vector(1,0){.07}}\put(52,70.25){\line(1,0){24.875}}
\put(64.88,30.81){\vector(1,0){.07}}\multiput(52.38,30.88)(6.25,-.0313){4}{\line(1,0){6.25}}
\put(46.88,49.75){\vector(1,0){.07}}\put(38,49.75){\line(1,0){8.875}}
\put(64.44,57.94){\vector(1,0){.07}}\multiput(52,58)(6.2188,-.0313){4}{\line(1,0){6.2188}}
\put(64.63,43){\vector(1,0){.07}}\put(52.13,43){\line(1,0){25}}
\put(64.25,65.88){\vector(1,0){.07}}\qbezier(52,58)(64.06,73.81)(76.88,57.88)
\put(64.25,34.75){\vector(1,0){.07}}\qbezier(52.13,42.75)(63.88,26.69)(77.13,42.88)
\put(19.63,55.63){\makebox(0,0)[cc]{$T$}}
\put(14.25,50.38){\makebox(0,0)[cc]{$\tau $}}
\put(57.38,50.63){\makebox(0,0)[cc]{$\sigma $}}
\put(65.5,49.5){\makebox(0,0)[cc]{$A$}}
\put(65.5,41){\makebox(0,0)[cc]{$A$}}
\put(65.13,37){\makebox(0,0)[cc]{$B$}}
\put(65.13,28){\makebox(0,0)[cc]{$B$}}
\put(64.75,63.6){\makebox(0,0)[cc]{$T$}}
\put(64.75,72.6){\makebox(0,0)[cc]{$T$}}
\put(58,60){\makebox(0,0)[cc]{$\Sigma _1$}}
\put(58,40.5){\makebox(0,0)[cc]{$\Sigma _2$}}
\put(59,78.38){\makebox(0,0)[cc]{$\Delta ^\prime $}}
\put(59,24){\makebox(0,0)[cc]{$\Delta ^{\prime\prime }$}}

\put(19.63,45){\makebox(0,0)[cc]{$B$}}
\put(16,38){\makebox(0,0)[cc]{$\Delta ^{\prime\prime }$}}
\put(16,63){\makebox(0,0)[cc]{$\Delta ^\prime $}}
\end{picture}

  \caption{Shortening $t$--bands.}\label{pic2}
\end{figure}
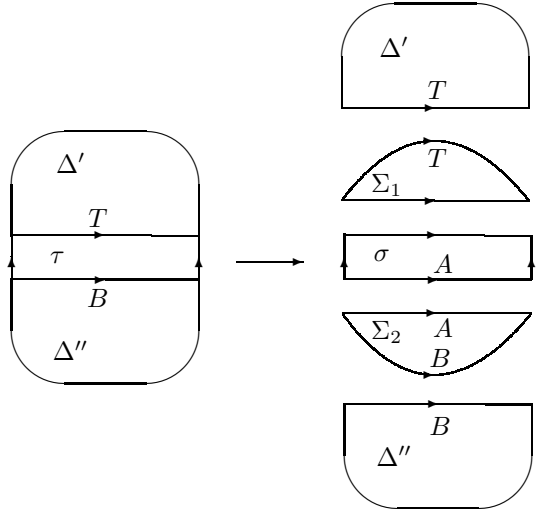

Let us prove the second assertion. Let $\Delta $ be an arbitrary
diagram over (\ref{p2}) with boundary label $W$, $\tau $  a
non--annular $t$--band in $\Delta $, $T\equiv \phi (\t (\tau ))$,
$B\equiv \phi (\b (\tau ))$. Let also $A$ be a shortest word in
$Y$ representing the same element of $K$ as $B$. Given these data,
we proceed as follows. First we remove $\tau $ from $\Delta $ thus
cutting the rest of $\Delta $ in two parts denoted by $\Delta
^\prime $ and $\Delta ^{\prime\prime }$. Further we take the
$t$--band $\sigma $ with $\phi (\b (\sigma ))\equiv A$. Note that
$$\phi (\t (\sigma ))= t^{-1}At=t^{-1}Bt=T$$ in $G$. Let $\Sigma
_1$, $\Sigma _2$ be diagrams over (\ref{p2}) having boundary
labels $\big( \phi (\t (\sigma ))\big) ^{-1}T$ and $AB^{-1}$
respectively. Gluing $\Delta ^\prime $, $\Sigma _1$, $\sigma $,
$\Sigma _2$, and $\Delta ^{\prime\prime }$ in the obvious way (see
Fig. \ref{pic2}), we obtain a new diagram with boundary label $W$,
where $$l(\sigma )=\| A\| =| k(\sigma )|_Y.$$ Obviously this
procedure does not violate condition 1) from the statement of the
lemma. Doing this for all non--annular $t$--bands, we get what we
need.
\end{proof}


\section{Proofs of the main results}


By technical reasons, it is convenient to introduce an auxiliary
presentation
\begin{equation}\label{p3}
G=\langle X\cup Y\cup \{ t\}\cup\mathcal H \; |\; \widetilde{
\mathcal R} \cup \mathcal Q \cup\mathcal S\rangle ,
\end{equation}
where $\mathcal Q$ is the set of all words in $Y$ representing $1$
in $K$, and $\mathcal S$, $\widetilde{ \mathcal R}$ are defined as
in (\ref{p1}) and (\ref{p2}). Speaking about van Kampen diagrams
over presentations (\ref{p1}), (\ref{p2}), and (\ref{p3}), we call
a cell $\Pi $ an $\mathcal R$--cell (respectively
$\widetilde{\mathcal R}$-- , $\mathcal R^\prime $--, $\mathcal
L$--, $\mathcal S$--, $\mathcal S_\lambda $--, $\mathcal Q$--cell)
if the boundary of $\Pi $ is labelled by a word from $\mathcal R$
(respectively $\widetilde{\mathcal R}$, $\mathcal R^\prime $,
$\mathcal L$, $\mathcal S$, $\mathcal S_\lambda $, $\mathcal Q$).
For a diagram $\Delta $, the set of all $\mathcal R$--,
$\widetilde{\mathcal R}$--, etc., cells (respectively the number
of such cells) is denoted by $\mathcal R(\Delta )$,
$\widetilde{\mathcal R} (\Delta )$, etc. (respectively $N_\mathcal
R(\Delta )$, $N_{\widetilde{\mathcal R}} (\Delta )$, etc.) Further
we set
$$M=\max\limits_{R\in \mathcal R^\prime } \| R\| .$$ Let
also $\gamma $ denote the relative Dehn function of $H$ with
respect to $\Hl\cup \{ K\} $ associated to (\ref{1}). In terms of
van Kampen diagrams this means that for any word $W$ in $X\cup
\mathcal H\cup (K\setminus \{ 1\} )$ representing $1$ in $H$,
there exists a diagram $\Delta $ over (\ref{p1}) such that
$N_{\mathcal R}(\Delta )\le \gamma (n)$.

\begin{lem}\label{WinXYH}
Let $W$ be a word in the alphabet $X\cup Y\cup \mathcal H$ of
length at most $n$ such that $W=1$ in $G$. Then there exists a van
Kampen diagram $\Delta $ over (\ref{p2}) with boundary label $W$
and the number of $\widetilde{\mathcal R}$--cells
\begin{equation} \label{Nrt}
N_{\widetilde{\mathcal R}} (\Delta ) \le n+(M+1)\gamma (n).
\end{equation}
\end{lem}

\begin{proof}
Note that we can regard $W$ as a word in $X\cup \mathcal H\cup
(K\setminus \{ 1\} )$ as $Y\subseteq K$. Obviously $W$ represents
$1$ in $H$. Let $\Delta _0$ be a diagram over (\ref{p1}) such
that:
\begin{enumerate}
\item[a)] $\partial \Delta _0\equiv W$;

\item[b)] The number of $\mathcal R$--cells in $\Delta _0$ is
$N_\mathcal R(\Delta _0)\le \gamma (n)$.

\item[c)] $\Delta _0$ has minimal total number of $\mathcal L$--
and $\mathcal S$--cells among all diagrams over (\ref{p1})
satisfying the first two conditions.
\end{enumerate}

\noindent Note that any internal edge of $\Delta _0$ (i.e., a
common edge of two cells) belong to the boundary of some $\mathcal
R$--cell. Indeed, if two $\mathcal S_\lambda $--cells $\Pi_1 $ and
$\Pi _2$ have a common edge $e$, we can erase $e$ replacing $\Pi
_1 $ and $\Pi _2$ with one cell since $\mathcal S_\lambda $
contains all words in $H_\lambda \setminus \{ 1\} $ representing
$1$ in $H$. However this contradicts c). The same argument can be
applied if two $\mathcal L $--cells of $\Delta _0$ have a common
edge.

For an edge $e$ in $\Delta _0$ labelled by a letter $k\in
K\setminus \{ 1\} $, we denote by $V(e)$ a shortest word in $Y$
representing $k$ in $G$. It is clear that replacing $e$ with
$V(e)$ for all such edges $e$ of $\Delta _0$, we obtain a diagram
$\Delta _1$ over (\ref{p3}). Note that $l(\partial \Delta
_1)=l(\partial \Delta _0)$ as any letter from $K\setminus \{ 1\} $
labelling an edge on $\partial \Delta _0$ belongs to $Y$. To
obtain a diagram over (\ref{p2}) it remains to get rid of
$\mathcal Q$--cells. For every $\mathcal Q$--cell $\Pi $ of
$\Delta _1$, we construct a diagram $\Sigma _\Pi $ over (\ref{p2})
as follows. Consider the $t$--annulus $\sigma $ whose outer
contour is $\b (\sigma )$ and $\phi (\b (\sigma ))\equiv \phi
(\partial \Pi )$. Obviously $\phi (\t (\sigma ))$ represents
$t^{-1}\phi (\Pi )t=1$ in $G$ and is a word in the alphabet $H_\nu
\setminus \{ 1\} $. Therefore, the we can glue an $\mathcal S_\nu
$--cell to the inner contour of $\sigma $. The resulting disk
diagram is denoted by $\Sigma _\Pi $.

Now replacing every $\mathcal Q$--cell in $\Delta _1$ with $\Sigma
_\Pi $, we obtain a diagram $\Delta $ over (\ref{p2}). Evidently
the number of $\widetilde{\mathcal R}$--cells of $\Delta $
satisfies
\begin{equation} \label{Nr}
N_{\widetilde{\mathcal R}} (\Delta )\le N_{\mathcal R} (\Delta
_0)+ \sum\limits_{\Pi \in \mathcal Q(\Delta _1)} l(\partial \Pi ).
\end{equation}
As any internal edge of $\Delta _0$ belongs to the boundary of
some $\mathcal R$--cell, any internal edge of $\Delta _1$ belongs
to the boundary of some $\mathcal R^\prime $--cell. Hence we have
\begin{equation}\label{slp}
\sum\limits_{\Pi \in \mathcal Q(\Delta _1)} l(\partial \Pi )\le
l(\partial \Delta _1)+\sum\limits_{\Pi \in \mathcal R^\prime
(\Delta _1)} l(\partial \Pi) \le \| W\| + MN_\mathcal R (\Delta
_0)\le n+M \gamma (n).
\end{equation}
Combining (\ref{Nr}) and (\ref{slp}) yields (\ref{Nrt}).
\end{proof}

\begin{proof}[Proof of Theorem \ref{HNN}]
Let $W$ be a word in the alphabet $X\cup Y\cup \{ t\} \cup
\mathcal H$ such that $\| W\| \le n $ and $W=1$ in $G$. We
consider a van Kampen diagram $\Delta $ that satisfies conditions
1), 2) from Lemma \ref{mainl}. To prove the theorem we have to
bound the number of $\widetilde{\mathcal R}$ cells in $\Delta $.

Let $D$ be a domain of $\Delta $. Since the boundary label of $D$
contains no letters $t^{\pm 1}$, we can think of $\phi (\partial
D)$ as a word in $X\cup\mathcal H\cup (K\setminus \{ 1\} )$.
Denote by $q_D$ an arbitrary cycle in $\G (H, X\cup\mathcal H\cup
(K\setminus \{ 1\}))$ having label $\phi (\partial D)$. For each
$\tau \in T(D)$, $\b (\tau )$ gives rise to a $K$--component of
$q_D$. Since no $t$--bands of $\Delta $ are $K$--connected, these
components are isolated in $q_D$.

Recall that for a $t$--band $\tau $, $k(\tau )$ denotes the
element of $K$ represented by $\phi (\b (\tau ))$. Let $c_D$
denote the cycle in $\G (H, X\cup\mathcal H\cup (K\setminus \{
1\}))$ obtained from $q_D$ by replacing all components with single
edges. We call the edges of $c_D$ corresponding to bottoms of
$t$--bands from $T(D)$ {\it distinguished}. Note that
distinguished edges are isolated $K$--components of $c_D$ and for
a $t$--band $\tau \in T(D)$, label of the distinguished edge
corresponding to $\tau $ represents the element $k(\tau )$.
Without loss of generality we may assume that $Y$ contains the set
$\Omega _K$ that is provided by Lemma \ref{Omega} applied to the
group $H$ and the collection of subgroups $\Hl \cup \{ K\} $.
Applying the second assertion of Lemma \ref{mainl} and Lemma
\ref{Omega}, we obtain
\begin{equation}\label{slt}
\sum\limits_{\tau \in T(D)} l(\tau ) = \sum\limits_{\tau \in T(D)}
|k(\tau )|_Y \le C\gamma (l(c_D)) .
\end{equation}

By $N_{top}(D)$ and $L_{\partial \Delta }(D)$ we denote the number
of $t$--bands $\tau $ of $\Delta $ such that $\t (\tau )\in
\partial D$ and the number of common edges of $\partial D$ and
$\partial \Delta $ respectively. Clearly $$l(c_D)\le
2(N_{top}(D)+L_{\partial \Delta }(D)).$$ Denote by $T(\Delta )$
the set of all non--annular $t$--bands in $\Delta $. Obviously we
have
\begin{equation}\label{D1}
\sum\limits_{D\in \mathcal D(\Delta )} l(c_D)\le \sum\limits_{D\in
\mathcal D(\Delta )}2(N_{top}(D)+L_{\partial \Delta }(D)) \le
2({\rm card\, } T(\Delta ) +n)\le 3n.
\end{equation}
Summing (\ref{slt}) over all domains of $\Delta $ and taking into
account (\ref{D1}), we obtain the following bound on the total
length of non--annular $t$--bands in $\Delta $
\begin{equation}\label{N1}
\sum\limits_{\tau\in T(\Delta )} l(\tau ) = \sum\limits_{D\in
\mathcal D(\Delta )} \sum\limits_{\tau \in T(D)} l(\tau ) \le
\sum\limits_{D\in \mathcal D(\Delta )} C\gamma (l(c_D))  \le
C\overline \gamma (3n).
\end{equation}

Furthermore, we have
\begin{equation}\label{ldD1}
\sum\limits_{D\in \mathcal D(\Delta )} l(\partial D ) \le
l(\partial \Delta ) + 2\sum\limits_{\tau\in T(\Delta )} l(\tau )
\le n+ 2C\overline \gamma (3n).
\end{equation}
By Lemma \ref{WinXYH} and inequality (\ref{ldD1}), we may assume
that the total number of all $\widetilde{\mathcal R}$--cells in
all domains of $\Delta $ is at most
\begin{equation}\label{N2}
\sum\limits_{D\in \mathcal D(\Delta )} \big[ l(\partial D) +
(M+1)\gamma (l(\partial D))\big] \le n +2C\overline \gamma (3n)
+(M+1)\overline\gamma (n+2C\overline \gamma (3n)).
\end{equation}
Finally summing (\ref{N1}) and (\ref{N2}) and taking into account
that $n\preceq f(n)$ for any superadditive function $f\not\equiv
0$, we obtain $N_{\widetilde{\mathcal R}} (\Delta )\preceq
\overline \gamma (n).$
\end{proof}

The next lemma is a relative analogue of a well--known property of
ordinary Dehn functions (see, for example, \cite{BMS}). The proof
is straightforward (and the same as in the non--relative case), so
we leave it to the reader.

\begin{lem}\label{retr}
Let $U$ be a group that is finitely presented relative to a
collection of subgroups  $\Ul $, and let $U_1$ be a retract of $U$
that contains all subgroups from the set $\Ul $. Suppose that the
relative Dehn function $\delta $ of $U$ with respect to  $\Ul $ is
well--defined. Then $U_1$ is finitely presented relative to $\Ul
$, the relative Dehn function $\delta _1$ of $U_1$ with respect to
$\Ul $ is well--defined and satisfies the inequality $\delta
_1\preceq \delta $.
\end{lem}

\begin{proof}[Proof of Theorem \ref{Am}]
Recall that the amalgamated product $A\ast _{K=\xi (K)} B$ is a
retract of the HNN--extension of the free product $A\ast B$ with
the associated subgroups $K$ and $\xi (K)$ \cite{LS}. It remains
to note that the relative Dehn function $\gamma $ of $A\ast B$
with respect to $\Am \cup \Bn $ satisfies (\ref{g}). This is
well--known for ordinary Dehn functions; the proof in the relative
case is obvious and actually the same, so we leave it as an
exercise to the reader. Now applying subsequently Theorem
\ref{HNN} and Lemma \ref{retr} we get Theorem \ref{Am}.
\end{proof}

\end{document}